\documentclass[12pt]{article}
\usepackage[utf8x]{inputenc}
\usepackage{amsmath}
\usepackage{amssymb}
\usepackage{graphicx}

\newcommand{\qed}{\hfill$\Box$\par\medskip\par\relax}

\numberwithin{equation}{section}

\newcommand{\eps}{\varepsilon}
\newcommand{\vr}{\varrho}
\newcommand{\Z}{{\mathbb Z}}

\newcommand{\Proj}{{\mathcal P}}
\newcommand{\LL}{{\mathcal L}}

\newcommand{\B}{{\mathcal B}}

\newcommand{\D}{{\mathcal D}}

\newcommand{\R}{{\mathbb R}}
\newcommand{\RR}{{\mathcal R}}

\newcommand{\X}{{\mathfrak X}}

\renewcommand{\phi}{\varphi}

\newcommand{\e}{\mathbf{e}}

\renewcommand{\S}{\mathcal{S}}
\newcommand{\Sph}{\mathbb{S}}

\newcommand{\IE}{\mathbb{E}}
\newcommand{\IP}{\mathbb{P}}

\allowdisplaybreaks

\newcommand{\F}{\mathcal{F}}

\newcommand{\1}[1]{{\mathbf{1}}_{\{#1\}}}

\newtheorem{theo}{Theorem}[section]
\newtheorem{lmm}[theo]{Lemma}
\newtheorem{df}[theo]{Definition}

\newtheorem{example}[theo]{Example}

\title{On range and local time of many-dimensional submartingales}
\author{Mikhail Menshikov$^{1}$ \and 
 Serguei~Popov$^{2,}$\footnote{Corresponding author}}

\begin{document}

\maketitle

{\footnotesize 
\noindent $^{~1}$University of Durham,
Department of Mathematical Sciences,
South Road, Durham DH1 3LE, UK.\\
\noindent e-mail:
\texttt{Mikhail.Menshikov@durham.ac.uk}

\noindent $^{~2}$Department of Statistics, Institute of Mathematics,
 Statistics and Scientific Computation, University of Campinas --
UNICAMP, rua S\'ergio Buarque de Holanda 651,
13083--859, Campinas SP, Brazil\\
\noindent e-mails: \texttt{popov@ime.unicamp.br}

}

\begin{abstract}
We consider a discrete-time process adapted to some filtration which lives
 on a (typically countable) subset of~$\R^d$, $d\geq 2$. For this process, we
assume that 
it has uniformly bounded jumps, is uniformly elliptic
 (can advance by at least some fixed amount with respect to any direction,
with uniformly positive probability). Also, we assume that the projection of
this process on some fixed vector is a submartingale, and that a
stronger additional
condition on the direction of the drift holds (this condition does not exclude that the drift could be equal to~$0$ or be arbitrarily small).
The main result is that with very high probability the number
of visits to any fixed site by time~$n$ is less 
than~$n^{\frac{1}{2}-\delta}$ for some $\delta>0$.
This in its turn implies that the number of different sites visited by
 the process by time~$n$ should be at least~$n^{\frac{1}{2}+\delta}$. 
\\[.3cm]\textbf{Keywords:} strongly directed submartingale,
Lyapunov function, exit probabilities
\\[.3cm]\textbf{AMS 2000 subject classifications:} 60G42, 60J10
\end{abstract}

\section{Introduction and results}
\label{s_intro}

Let $\X\subset\R^d$ be a set of infinite cardinality; the elements of~$\X$
will be 
called \emph{sites}. Without restriction of generality, we assume that
$0\in\X$. Throughout this
paper we assume that $d\geq 2$.  We consider a discrete-time process $X=(X_n,n\geq 0)$ with values in~$\X$,
adapted to a filtration $\F=(\F_n,n\geq 0)$. For the process~$X$,
we suppose that it is uniformly elliptic (can advance in any given 
direction with uniformly positive probability), has uniformly bounded
jumps, and is a \emph{strongly directed} submartingale (see Definition~\ref{d_terms} below for the precise meaning). In principle,
we do not assume homogeneity in space and/or time, or even the fact that the process is Markovian.

In this paper we study two related questions:
\begin{itemize}
 \item How many different sites can be visited by the process~$X$
by time~$n$?
 \item How large can be the number of visits to a given site?
\end{itemize}
Of course, in the absence of space/time homogeneity one cannot hope to be able to characterize the precise behavior of the quantities of interest; in this paper we content ourselves in proving that 
 with probability $1-\exp(-n^\eps)$ the number
of visits to any fixed site by time~$n$ is less 
than~$n^{\frac{1}{2}-\delta}$ for some $\delta>0$.
This in its turn implies that the number is different sites visited by the process by time~$n$ with very high probability will be at least~$n^{\frac{1}{2}+\delta}$.

Although it is not important for the formulation of our results,
while reading the paper one may always assume that~$\X$ is the vertex
set of the integer lattice~$\Z^d$, the vertex set of some other mosaic, or just any ``discrete'' (in particular, countable) set. This, of course, is justified by the questions that are of our interest: e.g., if the law of the jump of the process is (in some sense) continuous, then
such questions typically do not arise (every site is visited at most once, and the process visits~$n$ different sites by time~$n$). 

Range (i.e., the cardinality of the set visited sites, or sometimes 
this set itself) and the local time (i.e., the number of visits to a given
site) for space-homogeneous discrete-time random walks were extensively
studied in the literature. It is a classical result that the expected
 range of 
the simple random walk is $O(\frac{n}{\ln n})$ for $d=2$
and $O(n)$ for $d\geq 3$, see e.g.\ Section~6.1 of~\cite{Hughes}. 
It is not difficult to obtain from this fact (using an independence 
argument as e.g.\ in Lemma~3.1 of~\cite{AMP}) that
with very high 
probability the walk visits at least $n^{1-\delta}$ distinct sites
by time~$n$. 
Finer results for the range of homogeneous random walks can be found
in a number of papers; see e.g.\ in~\cite{BIK07,DV79,H98}
and references therein.
For nonhomogeneous random walks these questions, of course, are
more difficult; we mention~\cite{R07} 
that contains results on the range of
simple random walk on supercritical percolation cluster.  

The behaviour of the local time (i.e., the number of visits) in a fixed 
site, or the field of local times in all sites, was much studied 
in the literature as well. It is quite elementary to obtain that the expected number of visits to the origin by time~$n$ for the simple random walk is $O(\ln n)$ for $d=2$ and~$O(1)$ for $d\geq 3$.
Also, one can easily obtain for the simple random walk in dimension~$2$ (using e.g.~\textbf{E1} of Section~III.16 of~\cite{Spitzer}) that, with stretched-exponentially small probability,
the number of visits to the origin is less than~$n^\delta$ for any
fixed~$\delta>0$.
Of course, finer results (for more general random walks as well) 
are available; see e.g.\
\cite{Cerny,CCFR,CFR,CRR,MR}.

As mentioned above,
with our assumptions we cannot hope to obtain very ``precise'' results;
however, in some cases it may be important to be able to estimate the range from below. 
In particular, consider the following process, called 
\emph{excited random walk}, or sometimes \emph{cookie random walk}.
It is a discrete-time stochastic process  
taking values on~$\Z^d$, $d\geq 2$, described 
in the following way: when the particle
visits a site for the first time, 
it has a uniformly positive drift in
a given direction~$\ell$; when the particle is at a site which was
already visited before, it has zero drift
(observe that this implies that the cookie random walk is
a strongly directed submartingale in the sense of
Definition~\ref{d_terms} below). This process was introduced
in a simpler form ($\ell$~is the first coordinate vector and in already visited sites the process behaves as simple random walk)
in~\cite{BW03} and then studied (we mention only the papers that are concerned with dimension $d\geq 2$) in e.g.\ \cite{BR07,HH09,MPRV,Z06}.
The key fact that is usually needed is that such a process typically visits much more than~$n^{1/2}$ different sites by time~$n$;
this allows to prove bounds on the probability that the process advances
in the direction of cookies' drift by much more than~$Cn^{1/2}$, and 
this in turn makes is possible to use regeneration arguments that 
imply the Law of Large Numbers and the Central Limit Theorem. 

It is worth noting that in this paper, differently from what was
considered in~\cite{MPRV}, we allow the process to have a nonzero drift of arbitrarily small absolute value; for such processes, a direct application of methods of~\cite{MPRV} fails.

Now, we write formal definitions and state our results.
Let $\|\cdot\|$ be the Euclidean norm in $\R^d$ and let
$\Sph^{d-1}=\{x\in\R^d:\|x\|=1\}$ be the unit sphere.
Let us denote the coordinate vectors of~$\R^d$ by
$\e_1,\ldots,\e_d$.
We write $x\cdot y$ for the usual scalar product of $x,y\in\R^d$.
 For $A\subset \X$
we denote by~$|A|$ the cardinality of~$A$. Let
\begin{equation}
\label{d_localtime}
 L_n(x) = \sum_{k=0}^n \1{X_k=x}
\end{equation}
be the local time of the process in $x\in\X$ by time~$n$, 
and we denote by $\RR_n=\{X_0,\ldots,X_n\}$ the set of visited sites by time~$n$.
Define the random variable $D_n\in\R^d$ to be the (conditional) drift
of the process~$X$ at time~$n$:
\[
 D_n = \IE(X_{n+1}-X_n\mid \F_n) = \IE(X_{n+1}\mid \F_n)-X_n.
\]

\begin{figure}
 \centering \includegraphics{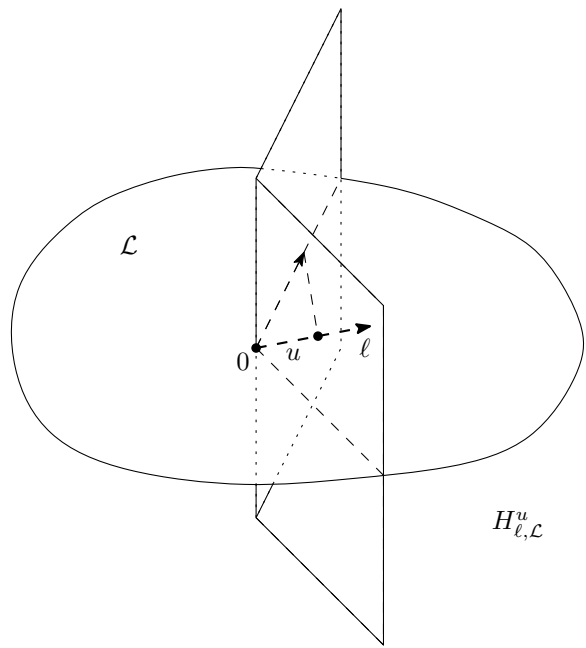} 
\caption{On the definition of $H_{\ell,\LL}^u$; observe that 
$H_{\ell,\LL}^u$ and $H_{-\ell,\LL}^{-u}$ ``complement'' each other (i.e., their union is~$\R^d$ and they intersect on a set of measure~$0$).}
\label{f_H}
\end{figure}

Next, let $\Proj_\LL$ be the operator of projection 
on the linear subspace $\LL\subset\R^d$.
Assuming that $\LL$ is a \emph{two-dimensional} subspace of~$\R^d$,
$\ell\in\Sph^{d-1}\cap\LL$ and $u\in\R$, define
\[
 H_{\ell,\LL}^u = \Big\{x\in\R^d: \Proj_\LL x=0 \text{ or }
         \frac{\Proj_\LL x\cdot\ell}{\|\Proj_\LL x\|}\geq u\Big\},
\]
see Figure~\ref{f_H}.

\begin{df}
\label{d_terms}
 We say that the $\F$-adapted process~$X$
\begin{itemize}
 \item[(a)] has uniformly bounded jumps, if there exists $K>0$ such that
$\|X_{n+1}-X_n\|\leq K$ a.s.\ for all~$n$ (we assume without restriction of generality that $K\geq 1$);
 \item[(b)] is uniformly elliptic (recall that we assume $d\geq 2$), if there exist $h,r>0$ such that
for all $\ell\in\Sph^{d-1}$ we have 
$\IP[(X_{n+1}-X_n)\cdot\ell > r\mid \F_n]>h$ a.s.\
(we assume without restriction of generality that $r\leq 1$);
 \item[(c)] is a martingale, if $D_n=0$ a.s.\ for all~$n$;
 \item[(d)] is $(u,\ell,\LL)$-strongly directed submartingale, if 
$u>0$ and $\IP[D_n \in H_{\ell,\LL}^u\mid \F_n]=1$ a.s.
\end{itemize}
\end{df}
Observe that for $(u,\ell,\LL)$-strongly directed submartingale 
it holds that the (conditional on the history) expected projection 
of the drift to~$\ell$ is always nonnegative (so that it is 
what one would naturally call a ``submartingale in direction~$\ell$''). 

In this paper we assume that all processes we are considering live
in a probability space with probability measure~$\IP$. We also adopt
the following notational convention: when it is necessary to assume that the initial state of the process under consideration is $x\in\X$,
we do not write it explicitly but simply add a subscript~$x$
to~$\IP$.

It is known that for the many-dimensional uniformly elliptic
martingales with bounded jumps the following
result (well hidden in~\cite{MPRV} as Lemma~5.3 and a part of
argument in Lemma~5.4) holds:
\begin{theo}
\label{t_mart}
Suppose that~$X$ is a uniformly elliptic martingale with uniformly bounded jumps.
Then, there exists $\hat\gamma\in (0,\frac{1}{2})$, 
${\hat C}_1,{\hat C}_2,
{\hat \delta}>0$ such that for any $x\in\X$
\begin{equation}
\label{eq_mart1}
 \IP_x[L_n(x)>n^{\hat \gamma}]
     \leq {\hat C}_1e^{-{\hat C}_2n^{\hat \delta}}
\end{equation}
and
\begin{equation}
\label{eq_mart2}
 \IP[|\RR_n|<n^{1-{\hat \gamma}}]
    \leq {\hat C}_1ne^{-{\hat C}_2n^{\hat \delta}}
\end{equation}
for all~$n$.
\end{theo}
One can note (we discuss this in detail in Section~\ref{s_proof_main}) that~\eqref{eq_mart2} follows from~\eqref{eq_mart1} in an
elementary way. 

The main purpose of this paper is to generalize Theorem~\ref{t_mart}
to (strongly directed) submartingales. Observe that one really needs some additional
condition on the submartingale; otherwise it may happen that
the typical number of visited sites is of order $n^{1/2}$. Indeed,
consider the following example:
\begin{example}
\label{ex_drift_to_0x}
 Let~$X$ be a nearest-neighbor random
walk on~$\Z^2$, with the transition probabilities described in the 
following way. From the horizontal axis, the particle
goes to neighboring sites with equal probabilities.
Off the horizontal axis, the particle always goes to the left/right with probabilities
$1/4$, the absolute value of the second coordinate increases with probability
$1/6$ and decreases with probability $1/3$.
This process is a submartingale in the direction of the 
first coordinate vector; on the other hand, it is clear (after some thought)
that the number of visited sites by time~$n$ behaves as 
$O(n^{1/2})$ (due to the drift towards the horizontal axis, the
process is essentially ``one-dimensional'', it ``lives'' in a small
neighborhood of the horizontal axis).
\end{example}

The main result of this paper is the following
\begin{theo}
\label{t_submart2}
 Suppose that $d\geq 2$ and~$X$ is a $(u,\ell,\LL)$-strongly directed submartingale, 
which is
 uniformly elliptic and has uniformly bounded jumps.
Then, there exist constants $\gamma\in (0,\frac{1}{2})$,
 $C_1,C_2,\delta>0$
(apart from the dimension, depending only on~$u$ and on $K,r,h$ from
  Definition~\ref{d_terms} (a)--(b)) such that for any $x\in\X$
\begin{equation}
\label{eq_submart1}
 \IP_x[L_n(x)>n^\gamma]\leq C_1e^{-C_2n^\delta}
\end{equation}
and
\begin{equation}
\label{eq_submart2}
 \IP[|\RR_n|<n^{1-\gamma}]\leq C_1ne^{-C_2n^\delta}
\end{equation}
for all~$n$.
\end{theo}


The key to the proof of Theorem~\ref{t_submart2} is a technical fact
about exit probabilities from two-dimensional rectangular domains
for martingales. 
For a simply connected domain~$\D\subset\R^2$, denote by
\[
 \tau(\D) = \min\{k\geq 1: X_k\notin\D\}
\]
the exit time from the domain~$\D$. Let 
\[
\S(\D) = \{\alpha X_{\tau(\D)-1}+(1-\alpha)X_{\tau(\D)}, 0\leq \alpha\leq 1\}\subset\R^2
\]
be the segment of the trajectory of the process at the instant 
when it leaves~$\D$.

\begin{figure}
 \centering \includegraphics{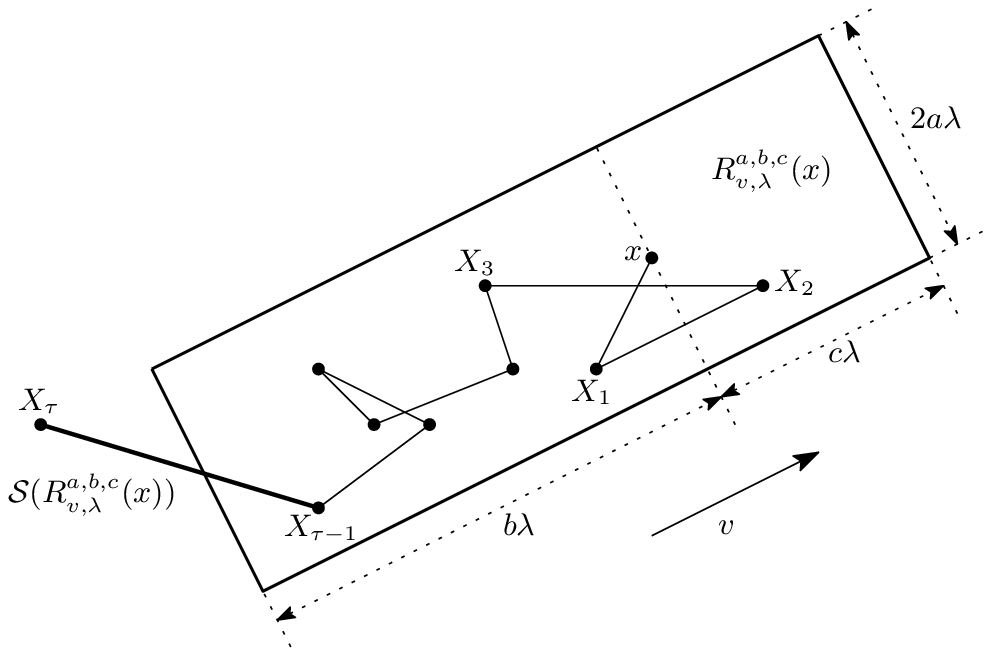} 
\caption{On the definition of the set $R_{v,\lambda}^{a,b,c}(x)$
and the event $G_{v,\lambda}^{a,b,c}(x)$, here we use the abbreviation 
$\tau:=\tau(R_{v,\lambda}^{a,b,c}(x))$.}
\label{f_rectangle}
\end{figure}

For $v\in\Sph^1$ and $a,b,c,\lambda>0$, let us define a rectangular domain
\[
 R_{v,\lambda}^{a,b,c}(x) 
        = \big\{y\in\R^2 : |(y-x)\cdot v^\bot|<a\lambda,
   (y-x)\cdot v\in(-b\lambda,c\lambda)\big\},
\]
where $v^\bot$ is a unit vector perpendicular to~$v$,
and consider the event
\[
 G_{v,\lambda}^{a,b,c}(x) = \Big\{\S(R_{v,\lambda}^{a,b,c}(x))
  \cap \{y\in\R^2 : |(y-x)\cdot v^\bot|\leq a\lambda,
   (y-x)\cdot v = -b\lambda\} \neq \emptyset\Big\}
\]
which means the process goes out of the rectangle through its 
``left'' side, see Figure~\ref{f_rectangle}.  

Now, the important technical fact in the proof of Theorem~\ref{t_submart2}
is
\begin{theo}
\label{t_exit_rectangle}
Assume that~$X$ is a uniformly elliptic martingale with uniformly bounded jumps.
 For all $a,b,c>0$ there exist $\vr=\vr(a,b,c)>0$ and $\lambda_0=\lambda_0(a,b,c)\geq 1$
such that 
\begin{equation}
\label{eq_exit_rectangle}
 \IP_x[G_{v,\lambda}^{a,b,c}(x)] \geq \vr
\end{equation}
for all $\lambda\geq \lambda_0$ and all~$x\in\X$.
\end{theo}
The explicit expressions for $\vr(a,b,c)$ and $\lambda_0(a,b,c)$
can be found at the end of the proof of Theorem~\ref{t_exit_rectangle},
see~\eqref{expl_lambda0} and~\eqref{expl_h}.

\begin{figure}
 \centering \includegraphics{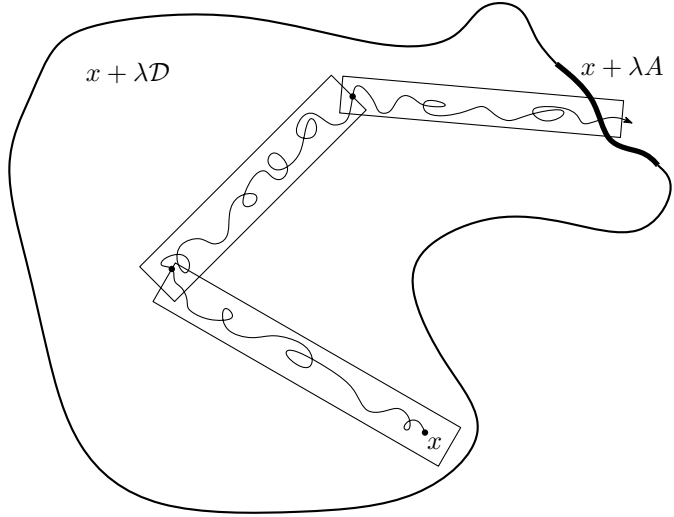} 
\caption{Going out from an arbitrary domain.}
\label{f_path_domain}
\end{figure}

In fact, there is only a small distance from the last result
to a more general one:
 let~$\D$ be a ``nice'' domain such
that $0\in\D$, and suppose that~$A\subset\partial\D$ is a connected
piece of the boundary of~$\D$ with positive measure. Then there
exist $\vr=\vr(\D,A)>0$ and $\lambda_0=\lambda_0(\D,A)$
such that 
$\IP_x[\S(x+\lambda\D)\cap (x+\lambda A)\neq \emptyset]\geq \vr$
for all $\lambda\geq \lambda_0$ and all~$x\in\X$.
It is not our intention in this paper to investigate, under which
(not very restrictive) precise geometric assumptions on~$\D$ and~$A$ the above fact holds; the way how it should follow from Theorem~\ref{t_exit_rectangle}
is (hopefully) made clear by Figure~\ref{f_path_domain}.

The paper is organized in the following way. In Section~\ref{s_exit_prob}
we prove Theorem~\ref{t_exit_rectangle}.
Then, in Section~\ref{s_proof_main} we prove
the main result of this paper, Theorem~\ref{t_submart2}.
In Section~\ref{s_final}, for comparison purposes we give a sketch of
the proof of Theorem~\ref{t_mart}, and then discuss some open problems.

\section{Exit probabilities from rectangles}
\label{s_exit_prob}
We begin with the following lemma.
\begin{lmm}
\label{l_thin_rectangle}
Let~$X$ be a uniformly elliptic martingale with bounded jumps, and
let $K,r,h$ be the corresponding constants (cf.\ Definition~\ref{d_terms} (a)--(b)). Assume that $b>0$ and 
$a\geq\frac{7K(b+K)}{r\sqrt{h}}$. Then, for all 
$x\in\X$ and $v\in\Sph^1$ we have
\begin{equation}
\label{eq_thin_rectangle}
 \IP_x[G^{a,b,b}_{v,\lambda}(x)] \geq \frac{1}{7}
\end{equation}
for all $\lambda\geq\frac{3K}{b}$.
\end{lmm}
\begin{figure}
 \centering \includegraphics{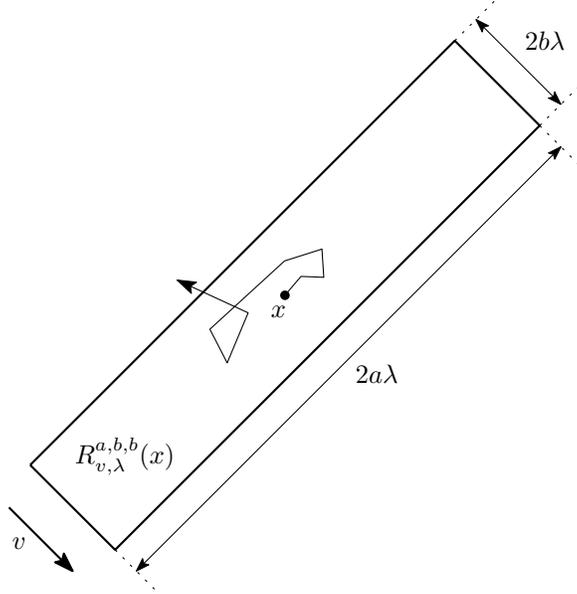} 
\caption{The event that the particle exits the
``stretched'' rectangle through its ``left'' long side.}
\label{f_exit_thin_rectangle}
\end{figure}
It is easy to believe that this result should hold true: it essentially means that
if the rectangle is (relatively) thin enough, then with uniformly
positive probability the process will exit the rectangle through a fixed long side (in this case, the ``left'' one), see Figure~\ref{f_exit_thin_rectangle}. In fact, the 
constant $1/7$ was chosen just for definiteness: one can modify
 the assumptions to obtain 
any fixed number less than $1/2$ on its place. 

\medskip
\noindent
\textit{Proof of Lemma~\ref{l_thin_rectangle}.}
Let us suppose without restriction of generality that $x=0$ and $v=\e_1$. Abbreviate $R^{a,b}_\lambda:=R^{a,b,b}_{\e_1,\lambda}(0)$.  Define
\begin{align*}
 \tau_\to &=\min\{k: X_k\cdot \e_1 > b\lambda\},\\
 \tau_\leftarrow&=\min\{k: X_k\cdot \e_1 < -b\lambda\},\\
 \tau_1 &=\min\{k: |X_k\cdot \e_1|> b\lambda\}
           =\tau_\to\wedge\tau_\leftarrow,
\end{align*}
and observe that $\tau_1\geq \tau(R^{a,b}_\lambda)$. 
By uniform ellipticity, it holds that 
$\IE((X_{m+1}\cdot\e_1-X_m\cdot\e_1)^2\mid\F_m)\geq r^2h$,
and this implies that the process $((X_m\cdot\e_1)^2-r^2hm,m\geq 0)$
is a submartingale with respect to the filtration~$\F$.
So, we have for any~$m$
\[
 \IE_0 \big((X_{m\wedge\tau_1}\cdot\e_1)^2-r^2h(m\wedge\tau_1)\big)
    \geq 0.
\]
Since, by the bounded convergence theorem,
\[
\lim_{m\to\infty} 
\IE_0(X_{m\wedge\tau_1}\cdot\e_1)^2 = \IE_0(X_{\tau_1}\cdot\e_1)^2 \leq (b\lambda+K)^2,
\]
and, by the monotone convergence theorem,
\[
\lim_{m\to\infty}  \IE_0(m\wedge\tau_1) = \IE_0\tau_1,
\]
we have $\IE_0\tau_1 \leq \frac{(b\lambda+K)^2}{r^2h}$. So, by
Chebyshev's inequality,
\begin{equation}
\label{tau1>}
\IP_0\Big[\tau_1 \geq \frac{7(b\lambda+K)^2}{r^2h}\Big] \leq \frac{1}{7}.
\end{equation}

Next, since $X_{\cdot\wedge\tau_1}\cdot\e_1$ is a (one-dimensional) 
bounded martingale, using the Optional Stopping Theorem we write
\begin{align*}
 0&= \IP_0[\tau_\leftarrow<\tau_\to] 
       \IE_0(X_{\tau_1}\cdot\e_1\mid \tau_\leftarrow<\tau_\to)
     + \IP_0[\tau_\leftarrow>\tau_\to] 
       \IE_0(X_{\tau_1}\cdot\e_1\mid \tau_\leftarrow>\tau_\to)\\
 &\geq -(b\lambda+K)\IP_0[\tau_\leftarrow<\tau_\to] 
         +b\lambda(1-\IP_0[\tau_\leftarrow<\tau_\to])\\
 &= b\lambda - (2b\lambda+K)\IP_0[\tau_\leftarrow<\tau_\to] ,
\end{align*}
so
\begin{equation}
\label{tau<tau}
 \IP_0[\tau_\leftarrow<\tau_\to] 
       \geq  \frac{b\lambda}{2b\lambda+K} > \frac{3}{7}
\end{equation}
since $\lambda>\frac{3K}{b}$.

Then, observe that Doob's inequality together
with the fact that the jumps are bounded by~$K$ imply that 
(abbreviate $s_\lambda:=\frac{7(b\lambda+K)^2}{r^2h}$)
\begin{align}
\IP_0\Big[\max_{j\leq s_\lambda}|X_j\cdot\e_2|\geq a\lambda\Big]
 &\leq \frac{\IE_0 (X_{\lfloor s_\lambda\rfloor}\cdot \e_2)^2}
{a^2\lambda^2} \nonumber\\
 &\leq \frac{K^2s_\lambda}{a^2\lambda^2} \nonumber\\
 &\leq \frac{7K^2(b+K)^2}{r^2ha^2} \leq \frac{1}{7},
\label{oc_vertical}
\end{align}
recall that we assumed that $a\geq \frac{7K(b+K)}{r\sqrt{h}}$.
The claim of Lemma~\ref{l_thin_rectangle} now follows from~\eqref{tau1>}, \eqref{tau<tau}, and~\eqref{oc_vertical}.
\qed

Now, we are ready to prove Theorem~\ref{t_exit_rectangle}.

\medskip
\noindent
\textit{Proof of Theorem~\ref{t_exit_rectangle}.}
For $x,y\in\R^2\setminus\{0\}$, let $\theta(x,y)\in [0,2\pi)$ be the angle
between~$x$ and~$y$ in the anticlockwise direction; for definiteness,
 we set $\theta(x,y)=0$ if at least one of the vectors $x,y$ equals~$0$. 
Let us abbreviate 
\begin{align}
\label{abbr1}
a_0 &:=\frac{7K(1+K)}{r\sqrt{h}},\\
\label{abbr2} 
\alpha_0 &:=\arctan\frac{1}{3a_0},\\
\label{abbr3}
m_0 &:=\Big\lceil\frac{\pi}{\alpha_0}\Big\rceil,\\
\label{abbr4}
 s_0 &:=2^{m_0}\times 6Ka_0 = {42}\frac{2^{m_0}K^2(1+K)}{r\sqrt{h}}.
\end{align}

Now, we define inductively two sequences of stopping
times: $\sigma(0)=T(0)=0$, and
\begin{align*}
 \sigma(m+1) &= \min\big\{k>\sigma(m): 
   \theta(X_{\sigma(m)},X_k)\in [\alpha_0,\pi)\big\}\\
 T(m+1) &= \min\Big\{k>\sigma(m): 
   X_k\cdot \frac{X_{\sigma(m)}}{\|X_{\sigma(m)}\|}
 \notin\Big[\frac{\|X_{\sigma(m)}\|}{2}, 
          \frac{3\|X_{\sigma(m)}\|}{2}\Big] \Big\},
\end{align*}
for $m\geq 1$.

\begin{figure}
 \centering \includegraphics{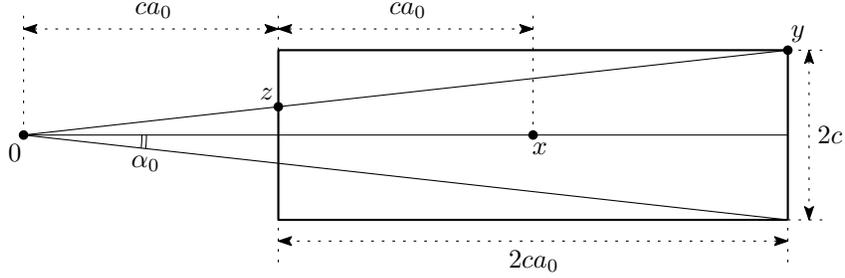} 
\caption{For $\alpha_0=\theta(x,y)$ we have $\alpha_0=\arctan\frac{1}{3a_0}$. Observe that $\|z\|\geq \|x\|/2$, $\|y\|\leq 2\|x\|$.}
\label{f_angle}
\end{figure}

Define the sequence of events $A_k=\{\sigma(k)<T(k)\}$, $k\geq 1$.
Lemma~\ref{l_thin_rectangle} then implies (see Figure~\ref{f_angle} and
recall that this lemma holds for rectangles with smaller side at least~$3K$)
that 
\begin{align}
 \IP\big[A_{k+1} \mid \F_{\sigma(k)}, \|X_{\sigma(k)}\|\geq 6Ka_0\big] 
   &\geq \IP\big[ G^{a_0,1,1}_{v_k,\lambda_k}(X_{\sigma(k)})
 \mid \F_{\sigma(k)}, \|X_{\sigma(k)}\|\geq 6Ka_0\big]\nonumber\\
           &\geq \frac{1}{7},\label{advance_angle}
\end{align}
where $v_k=\big(\frac{X_{\sigma(k)}}{\|X_{\sigma(k)}\|}\big)^\bot$,
$\lambda_k=\frac{\|X_{\sigma(k)}\|}{2a_0}$.

\begin{figure}
 \centering \includegraphics{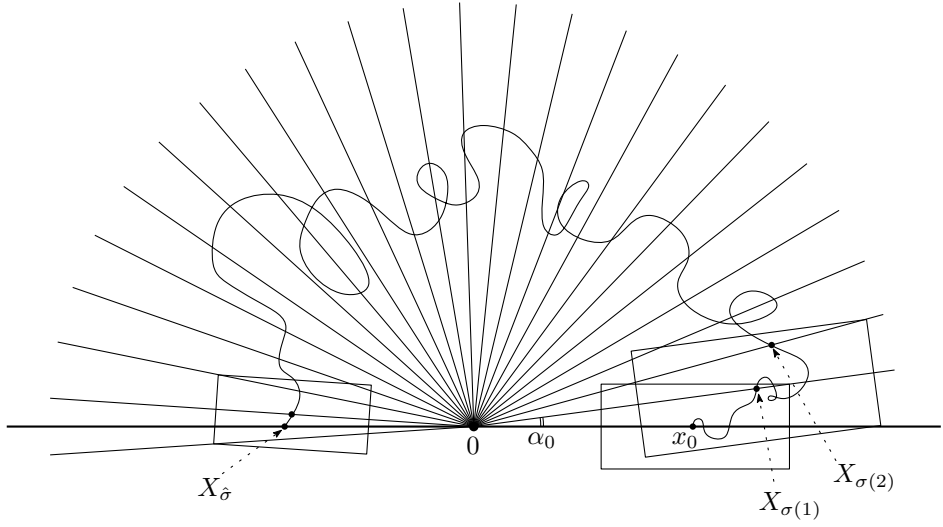} 
\caption{On the proof of~\eqref{eq_contornar}. For this picture, we
assume that the scale is such that~$K$ is not visible. The ratio 
of the longer side to the shorter side of the rectangles is~$a_0$.}
\label{f_contornar}
\end{figure}

Now, let us define 
\[
 {\hat\sigma} = \min\{k\geq 1 : X_k\cdot\e_1<0, 
        0\leq X_k\cdot\e_2< K\},
\]
and abbreviate as $\Proj_1$ the projector on the linear subspace spanned by~$\e_1$.
Then, \eqref{advance_angle} implies (see Figure~\ref{f_contornar})
that for any~$x_0$ such that $\|x_0\|\geq s_0$ (recall~\eqref{abbr4}),
$x_0\cdot\e_1>0$, $|x_0\cdot \e_2|\leq K$ we have
\begin{equation}
\label{eq_contornar}
 \IP_{x_0}\Big[X_{{\hat\sigma}}\cdot \e_1 \in
      [-2^{m_0}\|x_0\|,-2^{-m_0}\|x_0\|],\;
     \max_{k\leq{\hat \sigma}}\|X_k-\Proj_1 x_0\|<2^{m_0}\|x_0\|\Big]
  \geq \frac{1}{7^{m_0}}. 
\end{equation}

Now, we finish the proof of Theorem~\ref{t_exit_rectangle}.
Again, without lost of generality we assume that $x=0$ and $v=\e_1$.
Also, let us assume that $\lambda 2^{-m_0}(a \wedge c)>s_0$.

For $n\geq 1$, let us define the sequence of stopping times,
starting with $\tau_0=0$, by
\begin{align*}
 \tau_n &= \min\Big\{t \geq \tau_n :
 (X_{t}-X_{\tau_{n-1}})\cdot\e_1<
 -\lambda 2^{-m_0}(a \wedge c), \\
 &\qquad \qquad \qquad X_{t}\cdot\e_2\in[0,K],
  \max_{k\in [\tau_{n-1},t]}\|X_k-\Proj_1 X_{\tau_{n-1}}\|
     <\lambda(a \wedge c)\Big\}
\end{align*}
(we set formally $\min\emptyset=+\infty$),
and define the events $M_n=\{\tau_n<\infty\}$.
 From~\eqref{eq_contornar} we obtain that 
$\IP[M_n\mid \F_{\tau_{n-1}}]\geq 7^{-m_0}$, and, clearly, it holds that
\[
 G^{a,b,c}_{\e_1,\lambda}(0) \supset 
  \bigcup_{n=1}^{\lceil\frac{2^{m_0}b}{a \wedge c}\rceil} M_n
\]
(since on each of the events~$M_n$ the process advances to the left
by at least $\lambda 2^{-m_0}(a \wedge c)$, while staying inside
the rectangle~$R^{a,b,c}_{\e_1,\lambda}(0)$).
This concludes the proof of Theorem~\ref{t_exit_rectangle};
the explicit expressions for~$\vr$
and~$\lambda_0$ are then given by
\begin{align}
\label{expl_lambda0}
\lambda_0 &= \frac{s_0}{a \wedge c}\\ 
\intertext{and}
\label{expl_h}
 \vr &= \exp\Big(-m_0\Big\lceil\frac{2^{m_0}b}{a \wedge c}
                 \Big\rceil\ln 7\Big),
\end{align}
with $m_0,s_0$ defined by \eqref{abbr1}--\eqref{abbr4}.
\qed

\section{Proof of Theorem~\ref{t_submart2}}
\label{s_proof_main}

We consider first the case $d=2$. For this case, abbreviate
$H^u_\ell := H^u_{\ell,\R^2}$. Define $Y_0=X_0$, 
\[
 Y_n = X_n - \sum_{k=0}^{n-1} D_k
\]
for $n\geq 1$. Clearly, $Y$ is a martingale with jumps
uniformly bounded by $K'=2K$. To prove that~$Y$ is
uniformly elliptic, let us define
\[
 \tilde{D}_n = \begin{cases}
                 \frac{D_n}{\|D_n\|}, & \text{ on } \{D_n\neq 0\},\\
                  \e_1, & \text{ on } \{D_n=0\}.
               \end{cases}
\]
Then, observe that $Y_{n+1}-Y_n=X_{n+1}-X_n-D_n$,
so we have
\begin{equation}
\label{unif_ell_Y1}
 \IP[(Y_{n+1}-Y_n)\cdot Z\geq r \mid \F_n] = 
 \IP[(X_{n+1}-X_n)\cdot Z\geq r \mid \F_n] \geq h
\end{equation}
for $Z\in\{\tilde{D}_n^\bot,-\tilde{D}_n^\bot\}$, and
\begin{equation}
\label{unif_ell_Y2}
 \IP[(Y_{n+1}-Y_n)\cdot (-\tilde{D}_n)\geq r \mid \F_n] \geq 
 \IP[(X_{n+1}-X_n)\cdot (-\tilde{D}_n)\geq r \mid \F_n] \geq h.
\end{equation}
Using that~$Y$ is a martingale and~$\tilde{D}_n$ is $\F_n$-measurable,
we have $\IE\big((Y_{n+1}-Y_n)\cdot \tilde{D}_n\mid \F_n\big)=0$. Since~\eqref{unif_ell_Y2}
implies that $\IE\big((Y_{n+1}-Y_n)\cdot \tilde{D}_n\mid \F_n\big)^-
 = \IE\big((Y_{n+1}-Y_n)\cdot \tilde{D}_n\mid \F_n\big)^+ \geq rh$ 
and we have
also $\big((Y_{n+1}-Y_n)\cdot \tilde{D}_n\big)^+\leq 2K$, it holds 
that\footnote{it is elementary to obtain that for any random variable~$\xi$ with $0\leq \xi\leq a$ a.s.\
    and $\IE \xi\geq b$, it is true that $\IP[\xi\geq b/2]\geq b/(2a)$} 
\begin{equation}
\label{unif_ell_Y3}
 \IP[(Y_{n+1}-Y_n)\cdot \tilde{D}_n\geq rh/2] \geq \frac{rh}{4K}.
\end{equation}
Then, \eqref{unif_ell_Y1}--\eqref{unif_ell_Y3} imply that 
the process~$Y$ is uniformly elliptic with $r'=\frac{rh}{2\sqrt{2}}$ and
$h'=\frac{rh}{4K}$
(recall that we assume without 
restricting generality that $K\geq 1$).

Next, let us define
\begin{align*}
 {\tilde\sigma} &= \min\{j\geq 1: Y_j\in H_{-\ell}^u\},\\
 \sigma_0 &= \min\{j\geq 1: Y_j\cdot\ell \leq 0\},
\end{align*}
and, for $k\geq 1$
\[
 \sigma_k = \min\{j\geq 1: Y_j\cdot\ell \geq k\}.
\]

Now, our goal is to find a lower bound on the probability that, starting
from~$0\in\X$, the process avoids the origin for the next $k$~steps. 
Clearly, since $X_k-Y_k\in H_\ell^u$ a.s., for any~$k$ it holds that 
\begin{equation}
\label{X->Y}
 \IP_0[Y_1\notin H_{-\ell}^u,\ldots,Y_k\notin H_{-\ell}^u]
  \leq \IP_0[X_1\neq 0,\ldots, X_k\neq 0],
\end{equation}
so we concentrate on finding a lower bound for the term in the left-hand side of the above display.

For that, let us prove that there exists~$\beta>\frac{1}{2}$ and~$m_1$
such that for all $m\geq m_1$
\begin{equation}
\label{oc_beta}
 \IP_y[\sigma_{2m}<{\tilde\sigma}] \geq \beta
\end{equation}
for all~$y\in\X$ such that $y\cdot\ell \in [m,m+K)$.
First, note that the Optional Stopping Theorem implies analogously to~\eqref{tau<tau} that 
\begin{equation}
\label{gambler_1/2}
 \IP_y[\sigma_{2m}<\sigma_0]\geq \frac{m}{2m+K}
\end{equation}
for all~$y$ such that $y\cdot\ell \in[m,m+K)$.

\begin{figure}
 \centering \includegraphics{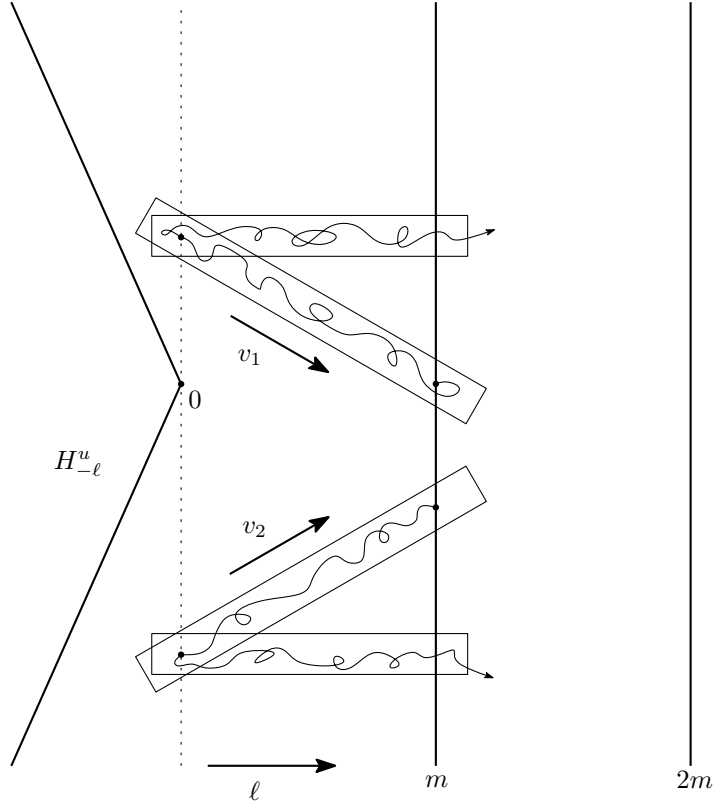} 
\caption{On the proof of~\eqref{oc_beta}}
\label{f_gamblers_ruin}
\end{figure}

Abbreviate $W_0^\ell = \{x\in\R^2: x\cdot\ell\in (-K,0]\}$ and
consider two vectors $v_1,v_2\in\Sph^1$ such that 
$v_j\cdot\ell>0$ for $j=1,2$ and $v_1+v_2$ is collinear to~$\ell$. Then, there exist $a_1,b_1,c_1,a_2,b_2,c_2>0$ such that 
\begin{itemize}
 \item[(i)] we have
\begin{align*}
 b_1|\cos \theta(v_1,\ell)|-a_1|\sin\theta(v_1,\ell)|&>1,\\
 a_1|\sin\theta(v_1,\ell)|+c_1|\cos \theta(v_1,\ell)|&<1
\end{align*}
 \item[(ii)] $b_2>1$,
 \item[(iii)] for $j=1,2$ the set 
\[R^{a_1,b_1,c_1}_{v_j,m}(m\ell) \cup
 \bigcup_{y\in W_0\cap R^{a_1,b_1,c_1}_{v_j,m}(m\ell)} 
     R^{a_2,b_2,c_2}_{-\ell,m}(y)\]
does not intersect $H_{-\ell}^u$,
\end{itemize}
see Figure~\ref{f_gamblers_ruin}.
In words, (i) means that (for large enough~$m$) the left side of the 
inclined rectangles
lies to the left of the dotted line (level~$0$) and the right side 
does not reach the level~$2m$, (ii) means that
the right side of the horizontal rectangle (relative to a point
 close to where the inclined rectangle intersect with the dotted line) is to
the right of the level~$m$, and (iii) means that the two rectangles cannot
touch the set $H_{-\ell}^u$. Then, it is clear that
Theorem~\ref{t_exit_rectangle} and~\eqref{gambler_1/2} imply that there exist
$m_1,\eps>0$
such that for all $m\geq m_1$
\begin{equation}
\label{go_and_be_back}
 \IP_y[\sigma_0<\sigma_{2m},\sigma_{2m}<{\tilde\sigma}] \geq \eps
\end{equation}
for any~$y$ such that $y\cdot\ell \in[m,m+K)$. 
Since
\[
 \IP_y[\sigma_{2m}<{\tilde\sigma}] = 
 \IP_y[\sigma_{2m}<\sigma_0]+ 
  \IP_y[\sigma_0<\sigma_{2m},\sigma_{2m}<{\tilde\sigma}],
\]
we obtain~\eqref{oc_beta} from~\eqref{gambler_1/2}
and~\eqref{go_and_be_back}.

Next, denote
\[
 k_1 = \Big\lceil\log_2\frac{K\sqrt{2n}}{m_1}\Big\rceil;
\]
observe that $m_1 2^{k_1}\geq K\sqrt{2n}$. Analogously to~\eqref{oc_vertical}, applying Doob's inequality we obtain
that
\begin{equation}
\label{2m->0}
 \IP_y[\sigma_0 > n] \geq \frac{1}{2}
\end{equation}
for all~$y$ such that $y\cdot \ell \geq K\sqrt{2n}$. Then,
\eqref{oc_beta} implies that for any~$y$ such that $y\cdot\ell\geq m_1$
\begin{equation}
\label{m1->Ksqrt2n}
 \IP_y[\sigma_{K\sqrt{2n}}< {\tilde\sigma}] \geq \beta^{k_1}.
\end{equation}
So, using \eqref{2m->0}--\eqref{m1->Ksqrt2n} and uniform ellipticity
(to assure that the process can initially advance to level~$m_1$),
we obtain for some ${\tilde c}>0$, $\eps>0$,
\begin{equation}
\label{not_return}
\IP_0[Y_1\notin H_{-\ell}^u,\ldots,Y_n\notin H_{-\ell}^u]
 = \IP_0[{\tilde\sigma}>n]
 \geq \frac{(h')^{\lceil m_1/r'\rceil}}{2} \beta^{k_1}
 \geq {\tilde c} n^{-\frac{1}{2}+\eps}
\end{equation}
since $\beta>\frac{1}{2}$. Using~\eqref{X->Y},
we obtain~\eqref{eq_submart1} for $d=2$; for $d\geq 3$ it then follows
if one considers the projection on~$\LL$. 

To prove~\eqref{eq_submart2}, define ${\hat\tau}_0=0$, 
\[
{\hat\tau}_{k+1} = \min\big\{m>{\hat\tau}_k: 
   X_m\notin\{X_0,\ldots,X_{m-1}\}\big\};
\]
 i.e., 
$({\hat\tau}_k, k\geq 0)$ is the sequence of times when the process
enters previously unvisited sites. Then, to obtain~\eqref{eq_submart2}, we use the fact that the process 
$X^{(k)}_\cdot = X_{\cdot+{\hat\tau}_k}$ satisfies conditions of the theorem, and apply the union bound
(again, the projection argument implies that~\eqref{eq_submart2} holds for all $d\geq 2$). 
This concludes the proof of Theorem~\ref{t_submart2}.
\qed

\section{Final remarks and open problems}
\label{s_final}

First, let us briefly sketch the proof of Theorem~\ref{t_mart}.
So, suppose that~$X$ is a martingale in dimension $d\geq 2$,
with bounded jumps and uniform ellipticity.
To begin, we show that there exist $b\in (0,1)$ 
close enough to~$1$ and $\gamma'>0$ 
(depending only on $K$, $h$, $r$ --- the constants in Definition~\ref{d_terms} (a)--(b)) such that
\begin{equation}
\label{Xb-submart}
\IE(\|X_{n+1}\|^b\mid \F_n) \geq \|X_n\|^b\1{\|Y_n\|>\gamma'}.
\end{equation}
To see that~\eqref{Xb-submart} holds, 
first observe that for a fixed~$y\in\R^d$ we have
\begin{align}
\label{norm_Taylor}
 \|x+y\|^b &= \big(\|x\|^2+2x\cdot y+\|y\|^2\big)^{b/2}\nonumber\\
&= 
\|x\|^b \Big(1+b\frac{x\cdot y}
{\|x\|^2}+\frac{b\|y\|^2}{2\|x\|^2}-\frac{1}{2}b(2-b)
\frac{(x\cdot y)^2}{\|x\|^4}+o(\|x\|^{-2})\Big), 
\end{align}
as $x\to \infty$. So, denoting by $\phi_n$ the angle between~$x$
and $\Delta_n:=X_{n+1}-x$, we have 
\begin{align}
\lefteqn{
 \IE(\|X_{n+1}\|^b - \|X_n\|^b \mid \F_n, X_n=x)
}\nonumber\\
  &= \frac{b}{2\|x\|^{2-b}} \big(\IE(\|\Delta_n\|^2\mid \F_n) 
    - (2-b)\IE(\|\Delta_n\|^2\cos^2\phi_n\mid \F_n)+o(\|x\|^{-2})\big).
\label{Taylor1}
\end{align}
Using the uniform ellipticity and the boundedness of jumps,
one can obtain that 
\[
 \IE(\|\Delta_n\|^2\cos^2\phi_n\mid \F_n) < (1-\eps')
      \IE(\|\Delta_n\|^2\mid \F_n)
\]
for some $\eps'>0$, so if $b<1$ is close enough to~$1$, 
the right-hand side of~\eqref{Taylor1} is positive for all
large enough~$x$ (see more details in the 
proof of Lemma~5.2 in~\cite{MPRV}).


Denote by $\B(x,s)=\{y\in \R^d: \|x-y\|\leq s\}$ the closed ball
of radius~$s$ centered in~$x$. The proof of Theorem~\ref{t_mart}
is now quite straightforward:
\begin{itemize}
 \item the Optional Stopping Theorem implies that, starting from~$x_0\in\R^d$, the process $X$ will
reach~$\R^d\setminus\B(x_0,\alpha)$ (without coming back to~$0$) with 
probability at least~$O(\alpha^{-b})$ (to apply the Optional Stopping Theorem, one has first to force the process a bit away from the origin, which happens with positive probability by uniform ellipticity).
\item Doob's inequality 
implies that from any place in~$\R^d\setminus\B(x_0,\alpha)$ with probability bounded away from~$0$ the process will not return to~$x_0$ 
after additional~${\tilde c}_1\alpha^2$ steps, where ${\tilde c}_1>0$ is a 
small enough constant.
\item Then, consider $\alpha={\tilde c}_2n^{1/2}$ 
for large enough~${\tilde c}_2$.
By the previous discussion, each time the process is in~$x_0$,
independently of the past it has probability at least
of order~$n^{-b/2}$
 of not returning to~$x_0$ during next~$n$ steps.
\item Take $\eps>0$ such that $b+\eps<1$.
Then, by an obvious coin-tossing argument,
the number of visits to~$x_0$ by time~$n$
will not exceed~$n^{\frac{b+\eps}{2}}$ with probability
at least $1-{\tilde c}_3e^{-{\tilde c}_4 n^{\eps/2}}$.
\item So, with probability
at least $1-{\tilde c}_3ne^{-{\tilde c}_4 n^{\eps/2}}$ 
the process~$X$ will have
 to visit at least $n^{1-\frac{b+\eps}{2}}$ different sites.
\end{itemize}

Now, we see that this proof is in sharp contrast with the proof of
Theorem~\ref{t_submart2}. Originally, our intention was to find 
a proof for Theorem~\ref{t_submart2} that would use Lyapunov functions
in a similar way; this amounts to finding a function~$f:\R^2 \mapsto\R$
with the following properties: $f\equiv 0$ on $\partial H^{-u}_{\e_1}$,
\[
\sup_{x\in \partial H^u_{\e_1}+n^{1/2}\e_1} f(x) = O(n^{\frac{1}{2}-\eps}) 
\]
for some $\eps>0$,
and $f(X_{\cdot\wedge\tau})$ is a submartingale, where $\tau$ is the hitting time of $H^u_{-\e_1}$. A possible idea would be to modify somehow
the function~$f_w$ of formula~(3.5) of~\cite{MMV}, but we did not 
succeed in developing it properly. 

\begin{figure}
 \centering \includegraphics{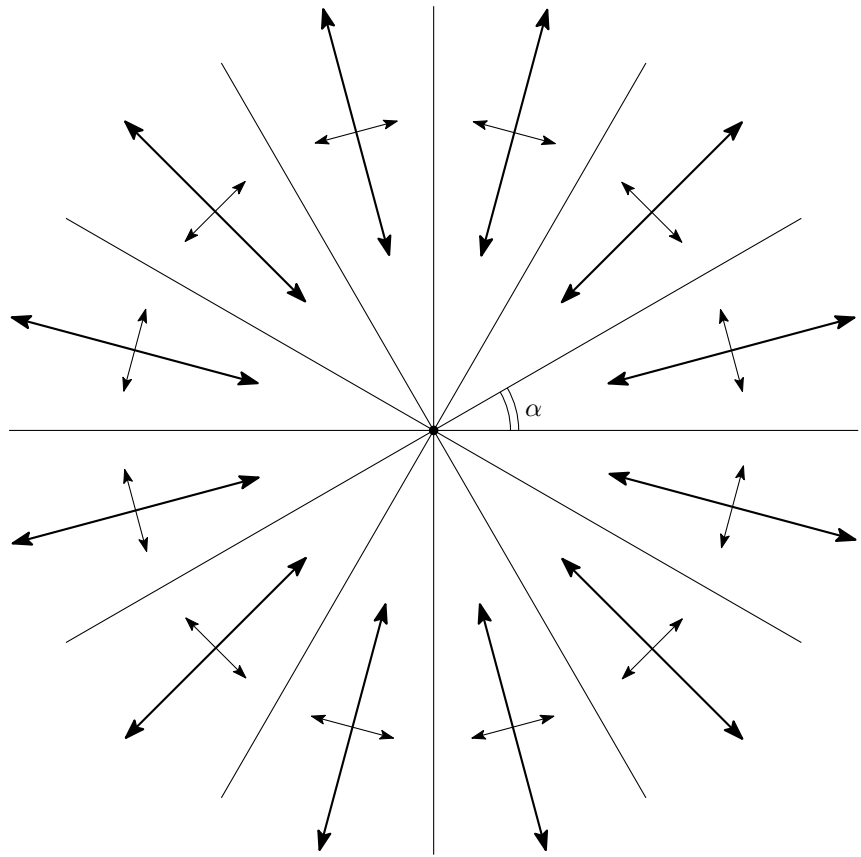} 
\caption{Example of a recurrent random walk with many visits to the origin.}
\label{f_recurrent_Zd}
\end{figure}

A natural question is if the results of Theorems~\ref{t_mart}
and~\ref{t_submart2} can be improved. As for the results for the 
local time (\eqref{eq_mart1} and~\eqref{eq_submart1}), it is not the case,
as the following example shows:
\begin{example}
\label{ex_manyvisits}
 Consider a zero-drift random walk on~$\Z^2$, defined in the way indicated on Figure~\ref{f_recurrent_Zd}. More specifically, 
we first divide the plane in sectors with (small enough) angle~$\alpha$, 
and then define the transition probabilities in each sector 
in such a way, roughly speaking, that the walk ``prefers'' the 
radial direction to the transversal one. For any fixed~$b<1$, it is clear that 
one can define the parameters of the model so that the 
expression in the parentheses in the right-hand side of~\eqref{Taylor1}
is negative for all large enough~$x$ 
(because the absolute value of the cosine in~\eqref{Taylor1} will 
typically be close to~$1$). Applying Theorem~1 of~\cite{AIM},
we obtain that $\IE_0\tau_0^{p}<\infty$ for any $p<\frac{b}{2}$,
where $\tau_0$ is the hitting time of the origin. But this implies that 
there typically will be~$n^p$ visits to the origin by time~$n$ (and that, with probability stretched-exponentially close to~$1$, the number of visits will be at least $n^{p-\delta}$, where~$\delta$ is an arbitrarily small positive number).
\end{example}

However, the situation with~\eqref{eq_mart2} and~\eqref{eq_submart2}
is less clear: it is an open problem to find out if~\eqref{eq_mart2} 
(respectively, \eqref{eq_submart2}) 
should be valid for all~$\hat\gamma>0$ 
(respectively, for all~$\gamma>0$). In fact, the authors were unable
to find any examples of uniformly elliptic martingales with uniformly 
bounded jumps for which the expected range is of order less than that
of the simple random walk (as mentioned in the introduction, it is
$O(\frac{n}{\ln n})$ for $d=2$
and $O(n)$ for $d\geq 3$).

\section*{Acknowledgements} The work of Serguei Popov was partially
supported by CNPq (300328/2005--2) and FAPESP (2009/52379--8).
 The work of Mikhail Menshikov was partially
supported by FAPESP (2011/07000--0). Also, we thank the anonymous referee 
for careful reading of the first version, which lead to many improvements.

\end{document}